\documentclass{article}
\usepackage{amssymb,amsfonts}
\newcommand\implies{\Rightarrow}
\newcommand{\includedin}{\subseteq}

\newcommand{\N}{\mathbb{N}}

\newcommand{\DD}{\mathcal{D}}
\newcommand{\MM}{\mathcal{M}}
\newcommand{\NN}{\mathcal{N}}

\def\newmathrm#1{\expandafter\newcommand\csname #1\endcsname{\mathrm{#1}}}

\newmathrm{GL}
\newmathrm{SL}
\newmathrm{PSL}
\newmathrm{Stab}
\newmathrm{Th}
\newmathrm{eq}
\newmathrm{rk}
\newmathrm{dcl}
\newmathrm{acl}
\newmathrm{ord}
\newmathrm{tp}
\newmathrm{stp}
\newmathrm{qftp}

\newcommand{\union}{\cup}

\let\intersect=\intersection

\newtheorem{lemma}{Lemma}[section]
\newtheorem{theorem}[lemma]{Theorem}
\newtheorem{prop}[lemma]{Proposition}

\newtheorem{fact}[lemma]{Fact}

\newtheorem{remark}[lemma]{Remark}

\newtheorem{claim}{Claim}[lemma]

\newtheorem{defn}[lemma]{Definition}
\newtheorem{definition}[lemma]{Definition}

\newcommand{\proof}{\smallskip\noindent {\it Proof\/}:\par}
\newcommand{\proofof}[1]{\smallskip\noindent {\it Proof of #1\/}:\par}
\newcommand{\qed}{\hfill{\vrule height 8pt depth .05pt width 5pt}\smallskip}

\begin{document}

\title{Borovik-Poizat rank and stability}
\author{Jeffrey Burdges\thanks{Supported by NSF Graduate Research Fellowship}\\
Department of Mathematics, Rutgers University\\
Hill Center, Piscataway, New Jersey 08854, U.S.A\\
e-mail: {\tt burdges@math.rutgers.edu}
\and
Gregory Cherlin\thanks{Supported by NSF Grant DMS 9803417}\\
Department of Mathematics, Rutgers University\\
Hill Center, Piscataway, New Jersey 08854, U.S.A\\
e-mail: {\tt cherlin@math.rutgers.edu}}
 
\maketitle

\section{Introduction}
 
Borovik proposed an axiomatic treatment of Morley rank in groups,
later modified by Poizat, who showed that in the context of groups the
resulting notion of rank provides a characterization of groups of
finite Morley rank \cite{Po}.  (This result makes use of ideas of
Lascar, which it encapsulates in a neat way.) These axioms form the
basis of the algebraic treatment of groups of finite Morley rank
undertaken in \cite{BN}. 

There are, however, ranked structures, i.e.~structures on which a
Borovik-Poizat rank function is defined, which are not
$\aleph_0$-stable \cite[p.~376]{BN}.  In \cite[p.~9]{Po} Poizat raised
the issue of the relationship between this notion of rank and
stability theory in the following terms: ``\dots un {\it groupe} de
Borovik est une structure stable, alors qu'un univers rang\'e n'a
aucune raison de l'\^etre \dots '' (emphasis added).  Nonetheless, we
will prove the following:
 
\begin{theorem}\label{superstable}
 A ranked structure is superstable.
\end{theorem}
 
An example of a non-$\aleph_0$-stable structure with Borovik-Poizat
rank 2 is given in \cite[p.~376]{BN}.  Furthermore, it appears that
this example can be modified in a straightforward way to give
$\aleph_0$-stable structures of Borovik-Poizat rank 2 in which the
Morley rank is any countable ordinal (which would refute a claim of
\cite[p.~373, proof of C.4]{BN}).  We have not checked the details. 
This does not leave much room for strenghthenings of our theorem.  On
the other hand, the proof of Theorem \ref{superstable} does give a
finite bound for the heights of certain trees of definable sets 
related to unsuperstability, as we will see in \S\ref{superstability}.

Since Shelah gave combinatorial criteria both for instability as well
as for unsuperstability in a stable context, to prove the theorem we
need only show that these criteria are incompatible with the
Borovik-Poizat rank axioms.  Now the rank axioms apply only to one
structure, while Shelah's criteria take their simplest form in a
saturated model.  There are two ways to bridge this gap.  Our first
proof worked directly within the model in which the rank function is
defined, paying attention in the process to the {\it uniformity} of
various first order definitions.  In the proof we give here, we first
extract the first order content of the rank axioms, then work with
them directly in a saturated model.

We will present the original rank axioms together with a few basic
consequences in \S2. 
Their first order content is analyzed in \S3; we
call the rank notions that result {\it $BP_0$-ranks}.  In \S4 we 
prove stability, and superstability is proved in \S5. 
The proof of superstability 
does not depend on the full strength of the axioms, so we will 
develop the basic facts about rank in a more general context adequate
for these applications.

For the stability and rank issues that concern us here, we may always
assume that the language of the structure involved is countable, and
we take advantage of this in \S3.

\section{The Borovik-Poizat Axioms and variations}
 
Let $\MM$ be a structure.
Let $\DD$ be the collection of parametrically  definable subsets of
$\MM^\eq$, i.e.~the sets and relations interpretable in $\MM$.
We say that $\MM$ is a {\it ranked structure} \cite[p.~57]{BN}
if there is a rank function
$\rk:(\DD - \{\emptyset\}) \to \N$
which satisfies the following axioms for all
$A,B \in \DD$. Such a rank will be called a $BP$-rank.
\begin{enumerate}
\item[Axiom $1$]  \textit{(Monotonicity of rank)}
  $\rk(A) \geq n + 1$ iff there are infinitely many pairwise disjoint,
  nonempty, definable subsets of $A$, each of rank at least $n$.
\item[Axiom $2$]  \textit{(Definability of rank)}
  If $f$ is a definable function from $A$ to $B$ then for each integer $n$
  the set $\{ b \in B : \rk (f^{-1}(b)) = n \}$ is definable.
\item[Axiom $3$]  \textit{(Additivity of rank)}
  If $f$ is a definable function from $A$ onto $B$, and for all $b\in B$
  we have $\rk(f^{-1}(b)) = n$, then $\rk(A) = \rk(B) + n$.
\item[Axiom $4$]  \textit{(Elimination of infinite quantifiers)}
  For any definable function $f$ from $A$ into $B$ there is an integer $m$ such
  that for any $b \in B$ the preimage $f^{-1}(b)$ is infinite whenever it
  contains at least $m$ elements.
\end{enumerate}
We adopt the convention that $\rk(\emptyset) = -\infty$.

These axioms are unnecessarily strong for our purposes.  
We prefer to work with the following weaker form of Axiom 1, and to
omit additivity of rank entirely:

\begin{enumerate}
\item[Axiom $1.1$] \textit{(Weak monotonicity)}
  $\rk(A\union B)=\max(\rk(A),\rk(B))$
\item[Axiom $1.2$] \textit{(Finite degree)}
  If there are infinitely many pairwise disjoint,
  nonempty, definable subsets of $A$, each of rank at least $n$, then
  $\rk(A) \geq n + 1$ 
\end{enumerate}

Let ``Axiom $1'$'' mean Axioms $1.1$ and $1.2$. 
A rank satisfying Axioms $1'$, $2$, and $4$ will be called a $BP'$-rank. 

Axiom $1.1$ is easily derived from full monotonicity as in
\textrm{\cite[Lemma 4.10 p.~59]{BN}}, and of course Axiom $1.2$ is
simply part of Axiom $1$, so a $BP$-rank is a $BP'$-rank. 

The following special case of Axiom $1.1$ is often used without comment:

\begin{fact}\label{rk-le} \textrm{\cite[Lemma 4.9 p.~59]{BN}}
In a ranked structure $\MM^\eq$, if $A \subseteq B$ are two definable sets
in $\MM^\eq$, then $\rk(A) \le \rk(B)$.
\end{fact}
 
Two of our axioms are more conveniently phrased in terms of {\it
uniformly definable} families of sets.  A family of sets $\{S_i:i\in
I\}$ is uniformly definable over the structure $\MM$ if there is a
single formula $\phi(x,y)$ defined in $\MM^\eq$ and a choice of
parameters $b_i$ ($i\in I$) such that $S_i=\phi[\MM^\eq,b_i]$ for all
$i$.  With the formula $\phi$ fixed, we write $S_b$ for
$\phi[\MM^\eq,b]$.  Note that the formula $\phi(x,y)$ involves no
parameters, as their places are filled by the variables $y$. 

\begin{prop}\label{definability}
Let $\MM$ be structure with a rank function satisfying Axioms $2$ and $4$,
let $D$ be a definable subset of $\MM^\eq$, and let $\{S_b\}_{b\in D}$ be
a uniformly definable family of sets over $\MM$. Then:
\begin{enumerate}
\item
(\textit{Definability of rank})
For each integer $n$,
$\{b \in D : \rk(S_b) = n\}$ is definable. \cite[Lemma 4.23 p.~66]{BN}
\item
(\textit{Elimination of infinite quantifiers})
There is a bound $m$, depending only on the formula $\phi$,
on the size of the finite members of the family
$\{S_b\}$.
\end{enumerate}
\end{prop}

A property of rank which is sometimes taken for granted 
is invariance under definable bijections.
We do not know if this holds for $BP'$-ranks in general, and we will not
assume it. 

\begin{defn} 
Let $\MM$ be a structure with a rank satisfying Axiom $1.1$, and let $A,B$
be two definable sets in $\MM^\eq$.
We write $A\equiv B$ if the rank of their symmetric difference
is less than the rank of their union. This may also be stated as
follows: they have the same rank $n$, 
and the rank of their symmetric difference is less than $n$.
\end{defn}

Observe that by Axiom $1.1$ this relation is an equivalence relation. 
 
\begin{remark}
In Axiom $2$, or equivalently in Proposition \ref{definability}, 
part (1), we refer to
definability with parameters from the model $\MM$. If the language
of $\MM$ is countable, then we may take the set of parameters involved
to be countable as well.
\end{remark}

\section{Canonical extensions of definable ranks}

Whenever one has a definable rank function on a structure $\MM$, it
has a canonical extension to any elementary extension of $\MM$;
details are given below.  We study the canonical extensions of
$BP$-ranks or $BP'$-ranks in the present section, giving axioms which
are satisfied by these ranks in general, and which exactly
characterize these canonical extensions in the case in which the
language is countable. 

The focus of interest is Axiom 1.2, which could be phrased as follows:
every set has a degree (analogous to the Morley degree).  We will
review the theory of the degree below, and show how to replace Axiom
1.2 by a degree approximation property which holds in canonical
extensions.  First we deal with the issue of the ``canonical
extension.'' 

\begin{defn}
Let $\MM$ be a structure equipped with a definable rank notion,
that is a function from definable subsets of $\MM^\eq$ to $\N$ such that
for any uniformly definable family $\{S_b\}$ over $\MM$, the set 
$$\{b:\rk(S_b)=n\}$$
is definable  (Axiom 2), 
and suppose that the rank has the following monotonicity
property: if $A\includedin B$ then $\rk(A)\le \rk(B)$ 
(a consequence of Axiom 1.1). This implies that the rank is bounded on
any uniformly definable family of sets.
Let $\NN$ be an elementary extension of $\MM$. 
The {\rm canonical extension} of $\rk$ to $\NN$ (again denote ``$\rk$'')
is defined as follows. Every parametrically definable subset $B$ of
$\MM^\eq$ has a canonical extension $B^*$ to $\NN^\eq$. 
For $S\includedin \NN^\eq$ definable, 
let $\rk(S)=n$ iff there is a uniformly definable
family $\{S_b:b\in B\}$ indexed by a parametrically definable subset
of $\MM$ such that $\rk(S_b)=n$ for $b\in B$, and $S=S_{b^*}$ for some
$b^*\in B^*$. (We have to check that this produces a well-defined rank
function.) 
\end{defn}

\begin{lemma}
With the hypotheses and notation of the preceding
definition, every definable subset of $\NN^\eq$ is assigned a
well-defined rank by the canonical extension of the rank function from
$\MM$ to $\NN$.
\end{lemma}
\proof
Let $S$ be definable in $\NN^\eq$. Then we can find a uniformly
definable family $(S_b:b\in B)$ with $B$ $0$-definable, such that
$S=S_{b^*}$ for some $b^*\in B$. Furthermore each set $S_b$ is
contained in a single sort of $\NN^\eq$, and by our (very weak) monotonicity
hypothesis, $\rk(S_b)$ is bounded for $b\in B$. Hence by definability
of ranks, $B^\MM$ can be partitioned into a finite number of
$\MM$-definable sets $B_i$
such that $\rk(S_b)=i$ on $B_i$. Then $b^*\in B_i^\NN$ for some $i$, and
hence we get $\rk(S_b)=i$. Thus every definable set is assigned at
least one rank. 

We must also verify that no conflicts arise. Suppose therefore that
$S$ lies in the extension to $\NN$ 
of the uniformly definable families $\{S_b:b\in B\}$ and $\{T_c:c\in
C\}$, and that $\rk$ is constant on both families. Since $\NN$ is an
elementary extension of $\MM$, it follows that in $\MM$ there are
some $b,c$ such that $S_b=T_c$. Hence the ranks are equal.
\qed

It is easy to see that any one of 
Axioms 1.1,2,3,4 will be preserved by the passage to a canonical
extension if it holds in the original model.
Our main interest at the moment will be in Axiom 1.2 (and later, in the
rest of Axiom 1).

\begin{defn}
If $\MM$ is a structure with a rank function satisfying Axiom 1.1, 
and $A$ is a definable set in
$\MM^\eq$, of rank $n$, then we say that $A$ has {\rm degree} $d$
if $A$ can be decomposed into $d$ disjoint definable pieces of rank
$n$, but no more. We say that $A$ has a degree (or, for emphasis, a
{\rm finite} degree) if it has degree $d$ for some $d$.  When $A$
has a degree, we use $\deg A$ to denote the degree of $A$.
\end{defn}
 
The theory of degree for the case of $BP$-rank is dealt with in 
\cite[Lemma 4.12,4.13]{BN}. In our context this reads as follows:

\begin{fact}\label{rk-equiv} 
Let $\MM$ be a structure with a rank function $\rk$ 
satisfying Axiom 1.1. Then
\begin{enumerate}
\item[1.] If $A$ is a definable set in $\MM^\eq$ which has rank $n$
and degree $d$,
then $A$ may be partitioned into $d$ definable pieces $A_i$ ($1\le
i\le d$) of rank $n$ and
degree $1$, and for any definable subset $B$ of $A$ of rank $n$ and
degree 1, we have $B\equiv A_i$ for a unique $i$.
In particular, the partition is unique modulo sets of lower rank 
in the following sense:
if $A$ is also decomposed as the union of $d$ definable subsets 
$A_i'$ of rank $n$ and degree 1, then after a permutation of the
indices, we will have $A_i\equiv A_i'$.
\item[2.] If $A,B$ are disjoint definable 
sets of equal rank which have degrees,
then $\deg(A\union B)=\deg(A)+\deg(B)$.
\item[3.] If the rank function satisfies Axiom $1'$ then every set has
a degree.
\end{enumerate}
\end{fact}

As a substitute for the existence of
degree, we consider the following, which can be stated loosely in the form:
``sets of finite degree are {\it dense}.'' 
As usual, we consider a rank function defined over a structure $\MM$;
and we also fix a set of parameters $C\includedin \MM$.

\begin{quote} \textit{(Degree Approximation Property over $C$)}
  If $A$ is a nonempty $C$-definable set, and $\{S_a:a\in A\}$ is a
  uniformly definable family (in particular, $S_a$ is $a$-definable for
  each $a$), then there is an element $a_0\in A$ for which $\deg(S_{a_0})$
  is finite.
\end{quote}

The case $C=\emptyset$ is reasonably strong but we will need this
relativized version. The case $C=M$ (the universe of $\MM$) is
pointless, as we are then assuming that every definable set has a
degree, since each definable set, taken by itself, constitutes a
$C$-definable family in this case.

\begin{defn}
A {\rm $BP_0'$-rank} on a structure $\MM$,
relative to a set of parameters $C$, is a function $f:\DD\setminus
\{\emptyset\}\to \N$ satisfying:
\begin{enumerate}
\item[Axiom 1.1] {\rm Weak Monotonicity}: $\rk(A\union B)=\max(\rk(A),\rk(B))$
\item[Axiom $1.2'$] {\rm The Degree $C$-approximation property}
\item[Axiom 2] {\rm Definability of rank, with parameters in $C$} 
(i.e.~Proposition \ref{definability}, with parameters in $C$).
\item[Axiom 4] {\rm Elimination of Infinite Quantifiers}
\end{enumerate}
\end{defn}

\begin{lemma}
Let $\MM$ be a structure equipped with a $BP'$-rank $\rk$, for which
the rank function is definable with parameters in $C\includedin \MM$.
Let $\NN$ be an
elementary extension of $\MM$. Then the canonical extension of $\rk$ to
$\NN$ is a $BP_0'$-rank relative to the same set of parameters.
Furthermore, the canonical extension satisfies the additivity axiom if
and only if the original rank function does.
\end{lemma}
\proof
In the context of a definable rank, the additivity property is clearly
first order, so the final claim is immediate.

We should however check the $C$-approximation property for degree.
As every nonempty $C$-definable set $B$ has a point in $\MM$, and
every $\MM$-definable set $S$ has a degree, we just have to check that
the degree of $S$ is unaltered by canonical extension. This reduces to
the case in which $S$ has degree $1$. So suppose that in $\NN$ we have
a definable subset $T$ of $S^*$ so that both $T$ and $S^*\setminus T$ have
rank $n=\rk(S)$. Using definability of rank we can pull this down to a
set $T_0$ defined in $\MM$ with $T_0$ and $S\setminus T_0$ of rank
$n$, a contradiction.
\qed

One may refine this slightly: 
the canonical extension of a $BP_0'$-rank is again a $BP_0'$-rank, by 
essentially the same argument. 
When the language is countable, our axioms actually characterize the
canonical extensions of $BP'$-ranks, as will be seen in the proof of
the next result. 

\begin{theorem}\label{ranks}
Let $T$ be a complete theory in a countable language. 
Then the following conditions are equivalent:
\begin{enumerate}
\item[$(A)$] $T$ has a model with a $BP'$-rank.
\item[$(B)$] $T$ has a model with a $BP_0'$-rank.
\item[$(C)$] $T$ has a countable extension by constants $T_C$, such that every
model of $T_C$ carries a $BP_0'$-rank, for which the degree approximation
property, and definability of rank, hold relative to the empty set.
\end{enumerate}
\end{theorem}

The equivalence of conditions $(B)$ and $(C)$ is clear, but worth noting,
and worth using: by passing to $T_C$ we can work over the empty set, and
lighten the notation. For the proof that $(C)$ implies $(A)$ we rely on the
following. 

\begin{lemma}\label{prime}
Let $T$ be a complete theory in a countable language, and suppose that $T$
has a model $\MM$ with a $BP_0'$-rank relative to $\emptyset$.
Then $T$ has a prime model.
\end{lemma}
\proof
We must show that $T$ is atomic, i.e. that every $\emptyset$-definable
nonempty set $X$ contains an $\emptyset$-definable atom.  
We may assume $X$ has minimal rank since ranks are finite.
By degree approximation $X$ has a degree,
since the family $\{X\}$ is already definable over the empty set.
We may suppose that $\deg(X)$ is minimized as well; we then claim that
$X$ is an atom over $\emptyset$.

By weak monotonicity any $\emptyset$-definable nonempty subset $Y$ of $X$
will have rank no larger than $\rk(X)$, and hence equal to $\rk(X)$ by the
minimization; accordingly  $Y$ will have finite degree no greater than
$\deg(X)$, and hence equal to $\deg(X)$; the same cannot apply
simultaneously to $X\setminus Y$, so $X\setminus Y$ must be empty, and
$Y=X$: $X$ is an atom over $\emptyset$.
\qed

\proofof {Theorem \ref{ranks}}
We assume $\NN$ is a structure on which we have a $BP_0'$-rank $\rk$
which has the degree approximation property and definability of rank
relative to the empty set. By Lemma~\ref{prime} the theory of $\NN$
has a prime model $\MM$, which we take to be an elementary
substructure of $\NN$. We claim that on $\MM$, the rank function gives
a $BP'$-rank. Only Axiom $1.2$ presents any issues: we claim that
every $\MM$-definable set $S$ has a degree.

Let $S$ be $a$-definable with $a\in M$ (an $n$-tuple for some $n$).
Let $A$ be the locus of $a$ over the empty set; as the type of $a$ is
principal, $A$ is a $0$-definable set. By the degree approximation
property, $A$ contains a point $a'$ for which the corresponding set
$S_{a'}$ has a degree $d$. As $\tp(a)=\tp(a')$ and rank is
$\emptyset$-definable, it follows easily that $\deg(S_a)=d$ as well.
\qed

We will show that the existence of $BP_0'$-rank implies
superstability. This is of course equivalent to the statement that
a $BP'$-rank gives superstability, since the problem localizes to
countable languages. As this is our main application, we have
emphasized $BP'$-ranks and $BP_0'$-ranks. However we can treat
$BP$-ranks similarly.

\begin{defn}
Let $\rk:\DD\setminus \{\emptyset\}\to \N$ 
be a rank function 
over a structure $\MM$, and $C\includedin \MM$ a set of parameters.
We say that $\rk$ has the {\rm splitting property} if every set of
rank $n>0$ contains a definable subset of rank $n-1$, and we say that 
$\rk$ has the {\rm splitting approximation property over $C$} if for
every uniformly definable family of infinite sets $\{S_b:b\in B\}$ 
indexed by a nonempty $C$-definable set $B$, there is an element $b\in B$
such that $S_b$ contains a definable subset $S'$ with 
$\rk(S')=\rk(S_b)-1$.
\end{defn}

Observe that a $BP'$-rank is a 
$BP$-rank if and only if it has the
splitting property and additivity. 
We define a $BP_0$-rank relative to a set of
parameters $C$ analogously, as a
rank satisfying Axiom $1.1$, degree and splitting approximation, 
and Axioms $2-4$, where the degree approximation property, the
splitting approximation property, and the
definability of rank all hold over $C$.

\begin{theorem}
Let $T$ be a complete theory in a countable language. 
Then the following conditions are equivalent:
\begin{enumerate}
\item[$(A)$] $T$ has a model with a $BP$-rank.
\item[$(B)$] $T$ has a model with a $BP_0$-rank.
\item[$(C)$] $T$ has a countable extension by constants $T_C$, such that every
model of $T_C$ carries a $BP_0$-rank, for which the degree approximation
property, definability of rank, and the splitting approximation
property all hold relative to the empty set.
\end{enumerate}
\end{theorem}
\proof
As before we need only prove $(C\implies A)$, and this reduces to the
claim that a $BP_0$-rank on the prime model is a $BP$-rank, with the
only property not yet verified being the splitting property.
Again, this reduces to the claim that if $\tp(a)=\tp(a')$ and
we have a uniformly definable family for which $S_{a'}$ contains a
definable subset $S'$ with $\rk(S')=\rk(S_{a'})-1$, then the same
applies to $S_a$. This is clear by definability of rank.
\qed

\section{Generic indistinguishability and Stability}\label{stability}
 
Before taking up stability as such, we analyze the structure of
definable binary relations in general. For this it will be convenient
to introduce a quantifier ``$\forall^* x$,'' read ``{\em for generic
$x$},'' as follows.

\begin{definition}
Assume $\MM$ carries a definable rank function satisfying Axiom 1.1.
Let $X$ be a definable set.
``$(\forall^* x \in X) \psi(x)$'' means: ``$\psi$ holds
generically on $X$,'' i.e.~$\rk(X - \psi[X]) < \rk(X)$.
By the definability of rank, if $X=X_a$ and $\psi=\psi(x,b)$ both vary
over uniformly definable families,
the set $\{(a,b):\hbox{$\forall^*x\in X_a$ $\psi(x,b)$}\}$
is definable. In other words, first order logic is closed under
the quantifier $\forall^*$.
\end{definition}

By Axiom 1.1,
if $\psi_1$ and $\psi_2$ hold generically on $X$,
then so does $\psi_1\&\psi_2$. On the other hand, the property:
``For all $\psi$, $\psi$ holds generically or $\neg \psi$ holds generically''
is equivalent to the condition that the degree of $X$ is equal to 1.
Note also that the relation $A\equiv B$ defined above can be
expressed as follows: $(\forall^* x\in A\union B) [x\in A\iff x\in B]$.
 
\begin{defn}
Let $\MM$ be a structure with a definable rank function satisfying
Axiom 1.1. Let $R$ be a definable relation on a definable set $S$.
We will say that $x_1,x_2 \in S$ are {\em generically indistinguishable}
for $R$ on $S$, and we write $x_1 \sim x_2$,
if $(\forall^* x \in S) (R(x_1,x) \iff R(x_2,x))$.
\end{defn}

Observe that this is an equivalence relation, and is definable
from whatever parameters are needed to define $R$ and $S$, 
together with those used to define rank.

\begin{prop}
Let $\MM$ be a structure with a $BP_0'$-rank.
Let $S$ be a definable set in $\MM^\eq$ and $R$ a definable binary
relation on $S$. Let $\sim$ be the relation of generic
indistinguishability for $R$ on $S$.
Then $S{/}{\sim}$ has finitely many classes.
\end{prop}

\proof
First, put $S$ and $R$ into a uniformly definable family 
$\{(S_a,R_a):a\in A\}$ 
with $A$ defined over the empty set. Let $\sim_a$ be the relation of
generic indistinguishability relative to
$R_a$ on $S_a$. Note that $\sim_a$ is definable from the parameter $a$
together with parameters needed to define certain ranks.
Then $\{S_a{/}{\sim}_a \}$ is a uniformly definable collection of sets.
Hence there is a uniform bound $m$ on the sizes of its finite members.
Consider  $I = \{a : |S_a{/}{\sim}_a| > m \}$, the set of indices for
which the quotient is infinite.  Our claim is that $I$ is empty.
If rank is $C$-definable, then $I$ is also $C$-definable.

Suppose $I$ is nonempty.
Then there is some $a \in I$ such that 
$d = \deg(S_a) < \infty$.  We will show that
$S_a{/}{\sim}_a$ has only finitely many classes to obtain a contradiction.

$S_a$ may be partitioned into $d$ definable pieces $S_{a,i}$ of degree $1$.
For $x_1\in S_a$, the set $\{x\in S_a:R_a(x_1,x)\}$
coincides with a union of some of the
$S_{a,i}$, modulo sets of lower rank. In other words, there is a set
$S'$, a union of finitely many of the $S_{a,i}$, for which:
$$(\forall^*x\in S_a) [R_a(x_1,x)\iff x\in S']$$
As there are at most $2^d$ possibilities for $S'$,
the relation $\sim_a$ has at most $2^d$ classes.
As this is finite, we have the desired contradiction.
\qed

\begin{theorem}
Let $\MM$ be a structure with a $BP_0'$-rank. Then $\Th(\MM)$ is stable.
\end{theorem}

\proof
We may replace $\MM$ by any elementarily equivalent model. So if the
theory of $\MM$ is unstable, we may suppose that there is a definable
relation $R$ on a definable subset $X$ of $\MM^\eq$ such that $R$
linearly orders some infinite subset of $X$, not necessarily definable.
We may also assume that $\MM$ is $\omega_1$-saturated.  Let $S$ be a
definable set which contains an infinite subset $L$ which is linearly
ordered by $R$, and has minimal rank.

Now we consider the relation $\sim$ of generic indistinguishability
for $R$ on $S$.  Since $S{/}{\sim}$ is finite, one of the equivalence
classes for $\sim$ on $S$ meets $L$ in an infinite set.  So without
loss of generality $S$ consists of a single $\sim$-class.  As $\MM$ is
$\omega_1$-saturated we may suppose $L$ has the order type of the
rationals.

Consider elements $a,b\in L$ with $a<b$.  The set $S'=\{x\in S:
[R(a,x) \iff \neg R(b,x)]\}$ contains the interval $(a,b)$ of $L$,
hence by the minimality of $\rk(S)$ we find
$$\rk(S')=\rk(S)$$
But this violates the generic indistinguishability of $a$ and $b$.
\qed

\section{Superstability}\label{superstability}
 
\begin{theorem}
Let $\MM$ be a structure with a $BP_0'$-rank. Then $\Th(\MM)$ is superstable.
\end{theorem}
 
\proof
In this proof we will be less cavalier about the distinction between
elements and $k$-tuples, as this will permit a slight refinement of the
result (see the remark following the proof).
 
Suppose $\MM$ is an unsuperstable structure with a $BP_0'$-rank.
As $\Th(\MM)$ is stable by the previous theorem, we can apply a
combinatorial criterion due to Shelah, involving an infinitely
branching tree of infinite height whose levels consist of pairwise
disjoint uniformly definable sets. This goes as follows.
 
In the first place we have
$D^1(x=x, L, \infty) = \infty$ \cite[Theorem II 3.14 p.~53]{Sh}.
Therefore, by
\cite[Lemma VII 3.5(5) p.~423]{Sh},
there are formulas $\phi_k \in L$ for $0\le k<\omega$
and in some model $\MM'$ of $T$ there are
parameters $a_v$ for $v$ a node of the tree $\omega^{< \omega}$,
such that:
 \begin{enumerate}
 \item[$(i)$] If $v\le w$ are two nodes in the tree, then
  $\phi_{|w|}[\MM',a_w] \includedin \phi_{|v|}[\MM',a_v]$ and
  $\phi_{|w|}[\MM',a_w] \ne \emptyset$;
 \item[$(ii)$]
   For any two distinct nodes $v,w$ at the same level $k$ of the tree,\\
       $\phi_k[\MM',a_v]\intersect \phi_k[\MM',a_w]= \emptyset$.
 \end{enumerate}
Note that the only condition imposed on the (single) root formula
$\phi_0[\MM',a_\emptyset]$
is that $\phi_0[\MM',a_\emptyset]$  should contain all
the sets $\phi_1 [\MM',a_{<i>}]$.
We may assume that $r_0 = \rk(\phi_0[\MM',a_\emptyset])$ is minimal among
all such trees.

For parameters $b$ of the same sort as $a_\emptyset$, consider the parameters
giving possible first level nodes of full rank
$$ N_b = \{c: \hbox{$\phi_1[\MM',c] \subseteq \phi_0[\MM',b]$
and $\rk(\phi_1[\MM',c])= r_0$} \} $$
Let $c_1 \sim_b c_2$ be the relation defined on $N_b$ by
$\phi_1[\MM',c_1] \equiv \phi_1[\MM',c_2]$.
$N_b$ and $\sim_b$ are uniformly definable
(by definability of rank) and $\sim_b$ is
an equivalence relation.

\begin{claim}
$N_a{/}{\sim}_a$ has finitely many classes
\end{claim}

\proof
First, $\{ N_x{/}{\sim}_x \}$ is a uniformly definable collection of sets,
so there is a uniform bound $m$ on the sizes of its finite members.  Let
$C = \{ x : |N_x{/}{\sim}_x| > m \}$ be the set of indices of infinite ones.
If $C \neq \emptyset$ then there is a $b \in C$ such that $d = \deg(\phi_0
[\MM',b]) < \infty$.  Now, $N_b{/}{\sim}_b$ has at most $2^d$ classes,
contradiction.
\qed

This indicates that there are only finitely many first level nodes of full
rank and there exists an index $i_0$ such that $\rk(\phi_1[\MM',a_{<i_0>}])
< r_0$.  We find a contradiction to the choice of $r_0$ by taking $\phi'_j
= \phi_{j+1}$ for $j\ge 0$ and $a_v' = a_{<i_0,v>}$.
\qed

\begin{remark}
We can prove a stronger result:
The height of any such tree of nonempty definable subsets of $\MM'$,
pairwise disjoint at each level, with infinite branching, cannot exceed
$r = \rk(\phi_0[\MM',a_\emptyset]) \leq \rk(\MM')$.  This goes by induction
on $r$ by following the line of the previous argument.
\end{remark}

\bibliographystyle{alpha}
 
\def\cprime{$'$}

\end{document}